\PassOptionsToPackage{unicode}{hyperref}
\PassOptionsToPackage{hyphens}{url}
\PassOptionsToPackage{dvipsnames,svgnames,x11names}{xcolor}
\documentclass[
  12pt]{article}

\usepackage{amsmath,amssymb, amsthm}

\makeatletter
\newcommand{\namelabel}[1]{%
  \phantomsection
  \renewcommand{\@currentlabel}{#1}
  \label{#1}
}
\makeatother

\usepackage{iftex}
\ifPDFTeX
  \usepackage[T1]{fontenc}
  \usepackage[utf8]{inputenc}
  \usepackage{textcomp} 
\else 
  \usepackage{unicode-math}
  \defaultfontfeatures{Scale=MatchLowercase}
  \defaultfontfeatures[\rmfamily]{Ligatures=TeX,Scale=1}
\fi
\usepackage{lmodern}
\ifPDFTeX\else  
\fi
\IfFileExists{upquote.sty}{\usepackage{upquote}}{}
\IfFileExists{microtype.sty}{
  \usepackage[]{microtype}
  \UseMicrotypeSet[protrusion]{basicmath} 
}{}
\makeatletter
\@ifundefined{KOMAClassName}{
  \IfFileExists{parskip.sty}{%
    \usepackage{parskip}
  }{
    \setlength{\parindent}{0pt}
    \setlength{\parskip}{6pt plus 2pt minus 1pt}}
}{
  \KOMAoptions{parskip=half}}
\makeatother
\usepackage{xcolor}
\makeatletter
\ifx\paragraph\undefined\else
  \let\oldparagraph\paragraph
  \renewcommand{\paragraph}{
    \@ifstar
      \xxxParagraphStar
      \xxxParagraphNoStar
  }
  \newcommand{\xxxParagraphStar}[1]{\oldparagraph*{#1}\mbox{}}
  \newcommand{\xxxParagraphNoStar}[1]{\oldparagraph{#1}\mbox{}}
\fi
\ifx\subparagraph\undefined\else
  \let\oldsubparagraph\subparagraph
  \renewcommand{\subparagraph}{
    \@ifstar
      \xxxSubParagraphStar
      \xxxSubParagraphNoStar
  }
  \newcommand{\xxxSubParagraphStar}[1]{\oldsubparagraph*{#1}\mbox{}}
  \newcommand{\xxxSubParagraphNoStar}[1]{\oldsubparagraph{#1}\mbox{}}
\fi
\makeatother

\usepackage{longtable,booktabs,array}
\usepackage{calc} 
\usepackage{etoolbox}
\makeatletter
\patchcmd\longtable{\par}{\if@noskipsec\mbox{}\fi\par}{}{}
\makeatother
\IfFileExists{footnotehyper.sty}{\usepackage{footnotehyper}}{\usepackage{footnote}}
\makesavenoteenv{longtable}
\usepackage{graphicx}
\makeatletter
\def\maxwidth{\ifdim\Gin@nat@width>\linewidth\linewidth\else\Gin@nat@width\fi}
\def\maxheight{\ifdim\Gin@nat@height>\textheight\textheight\else\Gin@nat@height\fi}
\makeatother
\setkeys{Gin}{width=\maxwidth,height=\maxheight,keepaspectratio}
\makeatletter
\def\fps@figure{htbp}
\makeatother

\makeatletter
\@ifpackageloaded{caption}{}{\usepackage{caption}}
\AtBeginDocument{%
\ifdefined\contentsname
  \renewcommand*\contentsname{Table of contents}
\else
  \newcommand\contentsname{Table of contents}
\fi
\ifdefined\listfigurename
  \renewcommand*\listfigurename{List of Figures}
\else
  \newcommand\listfigurename{List of Figures}
\fi
\ifdefined\listtablename
  \renewcommand*\listtablename{List of Tables}
\else
  \newcommand\listtablename{List of Tables}
\fi
\ifdefined\figurename
  \renewcommand*\figurename{Figure}
\else
  \newcommand\figurename{Figure}
\fi
\ifdefined\tablename
  \renewcommand*\tablename{Table}
\else
  \newcommand\tablename{Table}
\fi
}
\@ifpackageloaded{float}{}{\usepackage{float}}
\floatstyle{ruled}
\@ifundefined{c@chapter}{\newfloat{codelisting}{h}{lop}}{\newfloat{codelisting}{h}{lop}[chapter]}
\floatname{codelisting}{Listing}

\makeatother
\makeatletter
\@ifpackageloaded{caption}{}{\usepackage{caption}}
\@ifpackageloaded{subcaption}{}{\usepackage{subcaption}}
\makeatother

\ifLuaTeX
  \usepackage{selnolig}  
\fi
\usepackage[]{natbib}
\usepackage{bookmark}

\usepackage[normalem]{ulem}

\usepackage{pgfplots}
\usepackage{tikz}
\usepackage{tikz-3dplot} 
\usetikzlibrary{patterns,arrows, decorations.markings}
\usepgfplotslibrary{fillbetween}

\IfFileExists{xurl.sty}{\usepackage{xurl}}{} 
\hypersetup{
  pdftitle={A Duality of the F Distribution},
  pdfauthor={Zifeng Zhang; Peng Ding; Wen Zhou; Chenlu Shi},
  pdfkeywords={Analysis of variance; Fixed design; Linear regression; Random design; Spherical symmetry},
  colorlinks=true,
  linkcolor={blue},
  filecolor={Maroon},
  citecolor={Blue},
  urlcolor={Blue},
  pdfcreator={LaTeX via pandoc}}

\newcommand{\anon}{1}

\label{equation notation}
\numberwithin{equation}{section}
\def\be{\begin{equation}}
\def\ee{\end{equation}}
\def\bea{\begin{eqnarray}}
\def\eea{\end{eqnarray}}
\def\bd{\begin{displaymath}}
\def\ed{\end{displaymath}}
\def\bda{\begin{eqnarray*}}
\def\eda{\end{eqnarray*}}

\newcommand{\beq}{\begin{eqnarray*}}
\newcommand{\eeq}{\end{eqnarray*}}

\label{parenthesis notation}
\def\lsk{\left(}
\def\rsk{\right)}

\def\rbk{\right \} }
\def\lmk{\left [ }	
\def\rmk{\right ] }

\newcommand{\bone}{\mbox{\bf 1}}

\def\T{{ \mathrm{\scriptscriptstyle T} }}

\newcommand{\rank}{\mathrm{rank}}

\label{theorem notation}
\newtheorem{theorem}{Theorem}[section]
\newtheorem{lemma}{Lemma}[section]
\newtheorem{corollary}{Corollary}[section]

\newtheorem{proposition}{Proposition}[section]

\newcommand{\ie}{\mbox{\sl i.e.\;}}

\label{Var notation}

\newcommand{\p}{\mathbb{P}}

\label{ABC notation} 
\newcommand{\bA}{{\mathbf A}}

\newcommand{\bC}{{\mathbf C}}
\newcommand{\bD}{{\mathbf D}}

\newcommand{\bI}{{\mathbf I}}

\newcommand{\bP}{{\mathbf P}}
\newcommand{\bQ}{{\mathbf Q}}

\newcommand{\bW}{{\mathbf W}}
\newcommand{\bX}{{\mathbf X}}

\newcommand{\bZ}{{\mathbf Z}}

\label{abc notation}

\newcommand{\bfe}{{\mathbf e}}

\newcommand{\bu}{{\mathbf u}}

\newcommand{\bw}{{\mathbf w}}

\newcommand{\by}{{\mathbf y}}

\label{Greek letter notation}

\newcommand{\bbeta}  {\boldsymbol{\beta}}

\newcommand{\bepsilon}{\boldsymbol{\varepsilon}}

\newcommand{\bgamma}{\boldsymbol{\gamma}}
\newcommand{\btheta} {\boldsymbol{\theta}}
\newcommand{\bxi} {\boldsymbol{\xi}}

\newcommand{\bPhi} {\boldsymbol{\Phi}}

\newcommand{\bGamma} {\boldsymbol{\Gamma}}

\newcommand{\mt}{^{\T}}
\newcommand{\inverse}{^{-1}}

\newcommand{\bzero}{{\mathbf 0}}

\newcommand{\mbR}{\mathbb{R}}

\usepackage{accents}

\def\T{{ \mathrm{\scriptscriptstyle T} }}

\newcommand{\Perp}{\perp \! \! \! \perp}

\begin{document}

\def\spacingset#1{\renewcommand{\baselinestretch}%
{#1}\small\normalsize} \spacingset{1}


\if1\anon
{
  \title{\bf A Duality of the $F$ Distribution}
  \author{
  Zifeng Zhang\\
  \textit{Department of Biostatistics, Yale School of Public Health,}\\
  \textit{Yale University, New Haven, CT, USA}\\[0.75em]
  and\\[0.75em]
  Peng Ding\thanks{Corresponding author. E-mail: pengdingpku@berkeley.edu.}\\
  \textit{Department of Statistics, University of California, Berkeley,}\\
  \textit{Berkeley, CA, USA}\\[0.75em]
  and\\[0.75em]
  Wen Zhou\\
  \textit{Department of Biostatistics, School of Global Public Health,}\\
  \textit{New York University, New York, NY, USA}
  \\[0.75em]
  and\\[0.75em]
  Chenlu Shi\\
  \textit{Department of Mathematical Sciences,}\\
  \textit{New Jersey Institute of Technology, Newark, NJ, USA}
  }

  \date{October 19, 2026}
  \maketitle
} \fi

\if0\anon
{
  \bigskip
  \bigskip
  \bigskip
  \begin{center}
    {\LARGE\bf Title}
\end{center}
  \medskip
} \fi

\bigskip
\begin{abstract}
The $F$ distribution arises naturally in linear regression or analysis of variance. Elementary textbooks state the result in terms of the standard normal random vector and fixed projection matrices. More advanced textbooks and research articles have extended the result to spherically symmetric random vectors. We complement the literature with a dual result that the $F$ distribution arises with a fixed vector and random projection matrices generated by a left-spherical random matrix.  
\end{abstract}

\noindent%

\spacingset{1.8} 

\section{Overview of the Duality Result on the $F$ Distribution}
\label{sec: intro}

The $F$ distribution $F(v_1, v_2)$, due to \citet{Fisher24}, is defined as 
$$ 
F 
=
\frac{  \chi^2(v_1) / v_1 }{ \chi^2(v_2) / v_2 } ,
$$
where
$\chi^2 ( v_1 )$ and $\chi^2 ( v_2 )$ are independent $\chi^2$ distributions with degrees of freedom $v_1$ and $v_2$, respectively \citep{casella2024statistical}. 
It arises naturally in linear regression or analysis of variance. A classic result is stated below.

\begin{theorem}\label{thm::classic-F-1}
Suppose the following two assumptions hold. 
\begin{enumerate}
\item 
[(Z1)]\namelabel{(Z1)}
$\bZ$ is a fixed $n\times p$ matrix with $\mathrm{rank}(\bZ) = p$. 
\item 
[(W1)]\namelabel{(W1)}
$\bw$ is the $n$-dimensional standard normal random vector.
\end{enumerate}

Then we have
\begin{equation}
F 
=  
\frac{ \bw\mt (\bP_\bZ - \bP_{\bZ_2}) \bw / q}{ \bw\mt (\bI_n -\bP_\bZ) \bw / (n-p)} \sim F(q, n-p),
\label{eqn::F-normal}   
\end{equation}
where we partition $\bZ$ as $\bZ = [ \bZ_1\ \bZ_2 ]$ with $\bZ_1 \in \mathbb{R}^{n \times q}$ and $\bZ_2 \in \mathbb{R}^{n \times (p - q)}$, and $\bP_\bA = \bA (\bA\mt\bA)\inverse \bA\mt$ denotes the projection matrix onto the column space of $\bA$ with full column rank. 
\end{theorem}

Due to the ratio form of the $F$ in \eqref{eqn::F-normal}, \eqref{eqn::F-normal} holds if we multiply $\bw$ by any random variable $r$, independent of $\bw$, with no mass at 0. By \citet[][Theorem 2.5]{fang2018symmetric}, $r \bw$ is a spherically symmetric random vector. Therefore, Theorem \ref{thm::classic-F-1} also holds more generally for a spherically symmetric random vector, as stated below for completeness. 

\begin{theorem}\label{thm::classic-F-2}
\eqref{eqn::F-normal} holds if \ref{(W1)} is replaced by 
\begin{enumerate}
\item 
[(W1-1)]\namelabel{(W1-1)}
$\bw$ is a spherically symmetric random vector in the sense that  $\bw$ and $\bPhi \bw$ have the same distribution for any orthogonal matrix $\bPhi \in \mathbb{R}^{n \times n}$ \citep[][Definition 2.1]{fang2018symmetric}; and
\item
[(W1-2)]\namelabel{(W1-2)}
$\p \{ \bw = \bzero_n \} = 0$.
\end{enumerate}
\end{theorem}

Both Theorems \ref{thm::classic-F-1} and \ref{thm::classic-F-2} are known \citep{dempster1969elements}, although Theorem \ref{thm::classic-F-1} appears more often in elementary textbooks. 
Without the condition (W1-2), if we define $0/0=0$, then \eqref{eqn::F-normal} becomes a mixture of a point mass at 0 and the $F$ distribution.

We present a dual result of Theorem \ref{thm::classic-F-2} with fixed $\bw$ and a left-spherically symmetric matrix $\bZ$ in the sense that $\bZ $ and $ \bPhi \bZ$ have the same distribution for any orthogonal matrix $\bPhi \in \mathbb{R}^{n \times n}$
\citep{dawid1977spherical, kariya1981robustnessm}.

\begin{theorem}\label{thm::F-new}
\eqref{eqn::F-normal} holds if (Z1) and (W1)  are replaced by 
\begin{enumerate}
\item
[(Z2)]\namelabel{(Z2)}
$\bZ$ is a left-spherically symmetric random matrix with $\p \{ \mathrm{rank}(\bZ) = p \} = 1$; and 
\item 
[(W2)]\namelabel{(W2)}
$\bw \neq \bzero_n$ is a fixed vector.
\end{enumerate}
\end{theorem}

Related geometrical intuition appeared before \citep[e.g.,][]{dempster1969elements, efron1969student, eaton1970hotelling}. 
Geometrically, projection of a spherically symmetric random vector onto fixed column spaces is equivalent to projection of a fixed vector onto  column spaces generated by a left-spherically symmetric random matrix. We state the result rigorously in Theorem \ref{thm::F-new}.




We now turn to the regression implications of the previous results.
The two regimes considered in Theorems~\ref{thm::classic-F-2} and \ref{thm::F-new}, with random $\bw$ or random $\bZ$, naturally correspond to {\it fixed-design} or {\it random-design} linear regression, respectively.
In fixed-design linear regression, the design matrix is fixed whereas the error term is random \citep{lei2018asymptotics}. 
In random-design linear regression, the design matrix is random whereas the error term can be either fixed or random \citep{zhang2025random}. We present the linear regression analogues of Theorems \ref{thm::classic-F-2} and \ref{thm::F-new} below.

Consider the linear model
\begin{equation}
    \by = \bX \bbeta + \bepsilon,
\label{eqn: by = bX bbeta + bepsilon}
\end{equation}
where $\by \in \mathbb{R}^n$ is the response vector, $\bX \in \mathbb{R}^{n \times p}$ is the design matrix, $\bepsilon \in \mathbb{R}^n$ is the error vector, and $\bbeta \in \mathbb{R}^p$ is the parameter vector to be estimated. Assume $\mathrm{rank}(\bX) = p$. Define the ordinary least squares (OLS) estimator as
\begin{equation}
    \widehat{\bbeta} = (\bX\mt \bX)\inverse \bX\mt \by.
\label{eqn: OLS estimator}
\end{equation}
Given the parameter of interest $\bC \bbeta$, where $\bC \in \mathbb{R}^{q \times p}$ has $\mathrm{rank}(\bC) = q \leq p$, the following $F$ statistic is the foundation for testing hypotheses or constructing confidence regions \citep{ding2026linear}: 
\begin{equation}
F 
=
\frac{ \{ (\widehat{\bbeta}-\bbeta)\mt \bC\mt \lmk \bC (\bX\mt\bX)\inverse \bC\mt \rmk\inverse \bC ( \widehat{\bbeta}-\bbeta) \} / q }{  [ \by\mt \lsk \bI_n-\bP_{\bX} \rsk \by  ] / (n-p)}.
\label{eqn: F statistic in regression}
\end{equation}

The classic result states that the $F$ in \eqref{eqn: F statistic in regression} follows an $F$ distribution if the error terms are independent and homoskedastic normal.

\begin{theorem}
\label{thm::classic-F-ols}
The $F$ in \eqref{eqn: F statistic in regression} follows $F(q, n-p)$ if
\begin{enumerate}
\item[(X1)]\namelabel{(X1)}
$\bX$ is fixed; and
\item[(e1)]\namelabel{(e1)}
$\bepsilon \sim N(0, \sigma^2 \bI_n)$.
\end{enumerate}
\end{theorem}

Similar to Theorem \ref{thm::classic-F-2}, the $F$ in \eqref{eqn: F statistic in regression} follows an $F$ distribution if the error vector is spherically symmetric with zero point mass at 0 \citep{dawid1977spherical, kariya1981robustnessm, kariya1981robustness,  jensen1979linear}. 

\begin{theorem}
\label{thm::classic-F-ols-2}
The $F$ in \eqref{eqn: F statistic in regression} follows $F(q, n-p)$ if Condition \ref{(e1)} in Theorem \ref{thm::classic-F-ols} is replaced by
\begin{enumerate}
\item[(e1-1)]\namelabel{(e1-1)}
$\bepsilon$ is a spherically symmetric random vector; and
\item[(e1-2)]\namelabel{(e1-2)}
$\p \{ \bepsilon = \bzero_n \} = 0$.
\end{enumerate}
\end{theorem}

In parallel with Theorem \ref{thm::F-new}, we obtain the corresponding regression result: with fixed $\bepsilon$, the $F$ in \eqref{eqn: F statistic in regression} follows $F(q, n-p)$ if $\bX$ is a left-spherically symmetric random matrix. 

\begin{theorem}
\label{thm::F-ols-new}
The $F$ in \eqref{eqn: F statistic in regression} follows $F(q, n-p)$ if Conditions \ref{(X1)} and \ref{(e1)} in Theorem \ref{thm::classic-F-ols} are replaced by
\begin{enumerate}
\item[(X2)]\namelabel{(X2)}
$\bX$ is a left-spherically symmetric random matrix with $\p \{ \mathrm{rank}(\bX) = p  \} = 1$; and 
\item[(e2)]\namelabel{(e2)}
$\bepsilon \neq \bzero_n$ is a fixed vector. 
\end{enumerate}
\end{theorem}

Taking $\bC = \bfe_j\mt$ for $j \in \{1,\dots, p\}$, we can derive the dual result for the $t$ statistics from linear regression.

The next section presents the proofs of the aforementioned theorems.
\texttt{R} code for numerically confirming the new results in Theorems~\ref{thm::F-new} and \ref{thm::F-ols-new} is available at 
\href{https://github.com/JornZhang/F_Duality}{F Duality}.

\section{Proofs of Theorems}
\label{sec: tech and proofs}


The proofs of Theorems \ref{thm::classic-F-2} and \ref{thm::F-new} are elementary. We give their full proofs below for completeness. Our proofs leverage the following basic properties of spherically symmetric random vectors. 
Let $\mathrm{Unif}(S^{n-1})$ denote the uniform distribution on the unit sphere $S^{n-1} = \{ 
\bw \in \mathbb{R}^n: \| \bw \|_2 = 1 \}$.

\begin{lemma}
[Theorem 2.5 (i)(iii) in \cite{fang2018symmetric}]\label{prop: sufficient for Unif(S^(n-1))}
If $\bw \in \mbR^n$ is a spherically symmetric random vector, then we can represent it as $\bw \overset{d}{=} r \btheta$, where $r \Perp \btheta$, $\btheta \sim \mathrm{Unif}(S^{n-1})$,  and $r \in [0, \infty)$ is a random variable. 
\end{lemma}

In the special case with $\bxi \sim \mathcal{N} (\bzero_n, \bI_n)$, Lemma \ref{prop: sufficient for Unif(S^(n-1))} has a more transparent form: we can decompose it as $\bxi  = \| \bxi  \|_2 \cdot \bxi  / \| \bxi  \|_2 $, with  $\| \bxi  \|_2 \Perp \bxi  / \| \bxi  \|_2$, $\| \bxi  \|_2^2 \sim \chi^2 (n)$ and $\bxi  / \| \bxi  \|_2 \sim \mathrm{Unif}(S^{n-1})$.

Combining Lemma \ref{prop: sufficient for Unif(S^(n-1))} with the above representation for the standard normal random vector, we have the following result.

\begin{lemma}
\label{lemma::spherical-normal-representation}
If $\bw \in \mbR^n$ is a spherically symmetric random vector and $\bxi \sim \mathcal{N} (\bzero_n, \bI_n)$, then we can represent them via 
 $\bw \overset{d}{=} r \cdot \bxi / \| \bxi \|_2$, where $r \Perp \bxi$ and $\bxi \sim \mathcal{N} (\bzero_n, \bI_n)$.
\end{lemma}

We will use $\mathcal{C}(\bA)$ to denote the column space of $\bA$ and $\mathcal{C}(\bA)^{\perp}$ to denote its orthogonal complement. We will continue to use $\bP_\bA$ to denote the projection matrix onto $\mathcal{C}(\bA)$.

\subsection{Proof of Theorem \ref{thm::classic-F-2}}

It is a standard result but we still give a short proof for completeness.
By Lemma \ref{lemma::spherical-normal-representation}, we have $\bw \overset{d}{=} r \bxi / \| \bxi \|_2$, where $r \Perp \bxi$ and $\bxi \sim \mathcal{N} (\bzero_n, \bI_n)$.
The $F$ defined in \eqref{eqn::F-normal} satisfies
\begin{eqnarray}
F  
 \overset{d}{=} 
\frac{ \left. \left[ \frac{ r \bxi\mt}{\| \bxi \|_2} (\bP_\bZ - \bP_{\bZ_2}) \frac{ r\bxi}{\| \bxi \|_2} \right]  \middle/ q  \right.}{ \left. \left[ \frac{r \bxi\mt}{\| \bxi \|_2} (\bI_n -\bP_\bZ) \frac{r \bxi}{\| \bxi \|_2} \right] \middle/ (n-p) \right. } 
=
\frac{ \left[ \bxi\mt (\bP_\bZ - \bP_{\bZ_2}) \bxi \right]  / q}{ \left[ \bxi\mt (\bI_n -\bP_\bZ) \bxi \right] / (n-p)}.
\label{eqn: standard F with normal xi and fixed Z}
\end{eqnarray}

We can write $\bP_{\bZ} - \bP_{\bZ_2} = \bGamma_1 \bGamma_1\mt$ and $\bI_n - \bP_{\bZ} = \bGamma_2 \bGamma_2\mt$, where $\bGamma_1 \in \mathbb{R}^{n \times q}$ and $\bGamma_2 \in \mathbb{R}^{n \times (n-p)}$ are orthonormal bases of $\mathcal{C}(\bZ) \cap \mathcal{C}(\bZ_2)^{\perp}$ and $\mathcal{C}(\bZ)^{\perp}$, respectively, such that $\bGamma_1\mt \bGamma_2 = \mathbf{O}$.
Then we have $[\bGamma_1, \bGamma_2]\mt \bxi \sim \mathcal{N} (\bzero_{q+n-p}, \bI_{q+n-p})$.
By the definition of $\chi^2$ distribution and $F$ distribution, the $F$ in \eqref{eqn: standard F with normal xi and fixed Z} follows $F(q, n-p)$.

\subsection{Proof of Theorem \ref{thm::F-new}}

First, denote $\btheta = \bw / \| \bw \|_2$ so that $\btheta$ has $\| \btheta \|_2 = 1$. Then the $F$ defined in \eqref{eqn::F-normal} satisfies
\begin{eqnarray}
F 
= 
\frac{ \left. \left[ \frac{ \bw\mt}{\| \bw \|_2} (\bP_\bZ - \bP_{\bZ_2}) \frac{ \bw  }{\| \bw \|_2} \right]  \middle/ q  \right.}{ \left. \left[ \frac{ \bw\mt}{\| \bw \|_2} (\bI_n -\bP_\bZ) \frac{ \bw }{\| \bw \|_2} \right] \middle/ (n-p) \right. }  
=
\frac{ \btheta\mt (\bP_\bZ - \bP_{\bZ_2}) \btheta / q }{ \btheta\mt (\bI_n - \bP_\bZ) \btheta / (n-p) }. 
\label{eqn: F def with normalized vector}  
\end{eqnarray}
Although $\btheta \in \mathrm{S}^{n-1}$ under \ref{(W2)} is fixed, we can view it as a realization from $\boldsymbol{\vartheta} \sim \mathrm{Unif}(S^{n-1})$, which is independent of $\bZ$, i.e., $\boldsymbol{\vartheta} \sim \mathrm{Unif}(S^{n-1}) \Perp \bZ$.
Define
\begin{equation}
F' =   \frac{ \boldsymbol{\vartheta}\mt (\bP_\bZ - \bP_{\bZ_2}) \boldsymbol{\vartheta} / q }{ \boldsymbol{\vartheta}\mt (\bI_n - \bP_\bZ) \boldsymbol{\vartheta} / (n-p) } .
\label{eqn: F expression with randomized vartheta}
\end{equation}
By Theorem \ref{thm::classic-F-2}, we have $ F' \mid {\bZ} \sim F(q, n-p)$ for any given $\bZ$.
Hence we have $F' \Perp \bZ$ and the marginal distribution of $F'$ also follows $F(q, n-p)$.

Next, we show that $F' \Perp \boldsymbol{\vartheta}$ in \eqref{eqn: F expression with randomized vartheta} by proving $ \boldsymbol{\vartheta} \mid {F' = f'} \sim \mathrm{Unif}(S^{n-1})$ for any $f' \geq 0$.
For any orthogonal matrix $\bPhi \in \mathbb{R}^{n \times n}$, from \eqref{eqn: F expression with randomized vartheta}, we have
\begin{equation}
F'
= 
\frac{ \boldsymbol{\vartheta}\mt \bPhi \bPhi\mt(\bP_\bZ - \bP_{\bZ_2})\bPhi \bPhi\mt \boldsymbol{\vartheta} / q }{ \boldsymbol{\vartheta}\mt \bPhi \bPhi\mt(\bI_n - \bP_\bZ)\bPhi \bPhi\mt \boldsymbol{\vartheta} / (n-p) }  
=
 \frac{ \overline{\boldsymbol{\vartheta}}\mt (\bP_{\overline{\bZ}}-\bP_{\overline{\bZ}_2})\overline{\boldsymbol{\vartheta}} / q }{ \overline{\boldsymbol{\vartheta}}\mt (\bI_n - \bP_{\overline{\bZ}})\overline{\boldsymbol{\vartheta}} / (n-p) }  , 
 \label{eqn: F expression with arbitrary bPhi}
\end{equation}
where $\overline{\boldsymbol{\vartheta}}=\bPhi\mt\boldsymbol{\vartheta}$ and $\overline{\bZ} = \bPhi\mt \bZ$. 
By the spherical symmetry of $\mathrm{Unif}(S^{n-1})$ and the left-spherical symmetry of $\bZ$ under \ref{(Z2)}, we have $\overline{\boldsymbol{\vartheta}} \sim \mathrm{Unif}(S^{n-1})$, which is independent of  
$
\overline{\bZ}
=[
\bPhi\mt \bZ_1 \ \bPhi\mt\bZ_2
]
=
[
\overline{\bZ}_1 \ \overline{\bZ}_2
]
\overset{d}{=}
\bZ
$. 
Hence we have
\beq
F'
\overset{d}{=}
\frac{
\overline{\boldsymbol{\vartheta}}\mt
(\bP_{\bZ}-\bP_{\bZ_2})
\overline{\boldsymbol{\vartheta}}/q
}{
\overline{\boldsymbol{\vartheta}}\mt
(\bI_n-\bP_{\bZ})
\overline{\boldsymbol{\vartheta}}/(n-p)
}.
\eeq
Comparing this expression with
\eqref{eqn: F expression with randomized vartheta}, we have
\beq
(\boldsymbol{\vartheta} \mid F'=f')
\overset{d}{=}
(\overline{\boldsymbol{\vartheta}} \mid F'=f'),
\eeq 
that is, 
$$ 
(\boldsymbol{\vartheta} \mid {F' = f'} ) \overset{d}{=}  (\bPhi \mt \boldsymbol{\vartheta} \mid {F' = f'} )
$$ 
for any orthogonal matrix $\bPhi$. 
Thus $ \boldsymbol{\vartheta} \mid {F' = f'} $ is spherically symmetric with $\| \boldsymbol{\vartheta} \|_2 = 1$.
By Lemma \ref{prop: sufficient for Unif(S^(n-1))},  $ \boldsymbol{\vartheta} \mid {F' = f'} \sim \mathrm{Unif}(S^{n-1})$ for any $f' \geq 0$. 
Hence $F' \Perp \boldsymbol{\vartheta}$.

Finally, by $F' \Perp \boldsymbol{\vartheta}$, the conditional distribution $F'\mid {\boldsymbol{\vartheta} = \btheta} $ must be the same as  the marginal distribution $F(q, n-p)$, which is our desired result.

%

\subsection{Proof of Theorems \ref{thm::classic-F-ols-2} and \ref{thm::F-ols-new}}

Before we start the proof of Theorems \ref{thm::classic-F-ols-2} and \ref{thm::F-ols-new}, 
we first collect some results about the Frisch–Waugh–Lovell Theorem \citep{frisch1933part,lovell1963seasonal} summarized in \citet[Chap.~7]{ding2026linear}, which will be used frequently.

\begin{proposition}
[The Frisch–Waugh–Lovell Theorem]
Suppose that we have $\by \in \mathbb{R}^n$ and $\bW = [\bW_1, \bW_2] \in \mathbb{R}^{n \times p}$ with $\bW_1 \in \mathbb{R}^{n \times q}$, $\bW_2 \in \mathbb{R}^{n \times (p-q)}$, $\mathrm{rank}(\bW) = p$, and $1 \leq q \leq p \leq n$.
Define $\widehat{\bgamma} = (\bW\mt \bW)\inverse \bW\mt \by \in \mathbb{R}^p$ as the projection coefficient of $\by$ onto $\mathcal{C}(\bW)$, although there is not necessarily a linear model between $\bW$ and $\by$.
Also, $\widehat{\bgamma}$ has sub-vectors $\widehat{\bgamma} = [\widehat{\bgamma}_1\mt, \widehat{\bgamma}_2\mt]\mt$, where $\widehat{\bgamma}_1 \in \mathbb{R}^q$ and $\widehat{\bgamma}_2 \in \mathbb{R}^{p-q}$ correspond to $\bW_1$ and $\bW_2$, respectively.
Define $\widetilde{\by} = (\bI_n - \bP_{\bW_2}) \by$ and $\widetilde{\bW}_1 = (\bI_n - \bP_{\bW_2}) \bW_1$.
Then we have
\begin{enumerate}
    \item[(i)] $\mathrm{rank}(\widetilde{\bW}_1) = q$ and 
    $
    \widehat{\bgamma}_1 
    = 
    (\widetilde{\bW}_1\mt \widetilde{\bW}_1)^{-1} \widetilde{\bW}_1\mt \by 
    = 
    (\widetilde{\bW}_1\mt \widetilde{\bW}_1)\inverse \widetilde{\bW}_1\mt \widetilde{\by}
    $;

    \item[(ii)] For $(\bW\mt \bW)\inverse$, the top left $q \times q$ sub-matrix is $(\widetilde{\bW}_1\mt \widetilde{\bW}_1)\inverse$ and the bottom right $(p-q) \times (p-q)$ sub-matrix is $[\bW_2\mt (\bI_n - \bP_{\bW_1}) \bW_2]\inverse$;

    \item[(iii)] $\bP_{\bW} = \bP_{\bW_2} + \bP_{\widetilde{\bW}_1}$, $\bP_{\bW_2} \bP_{\widetilde{\bW}_1}  = \mathbf{O}$;

    \item[(iv)] $(\bI_n - \bP_{\bW}) \by = (\bI_n - \bP_{\widetilde{\bW}_1}) \widetilde{\by}$.
\end{enumerate}
\label{prop: FWL Thm}
\end{proposition}

Now we start the proof of Theorems \ref{thm::classic-F-ols-2} and \ref{thm::F-ols-new}.

Since $\rank(\bC)=q$, we can choose $\bD\in\mbR^{(p-q)\times p}$ such that $[\bC\mt \ \bD\mt]\mt \in \mathbb{R}^{p \times p}$ is invertible. 
Define
\begin{equation}
\bZ
=
\bX
\begin{bmatrix}
\bC\\
\bD
\end{bmatrix}^{-1}
=
[\bZ_1 \ \bZ_2],
\label{eqn: transformed Z in regression}
\end{equation}
where $\bZ_1\in\mbR^{n\times q}$ and $\bZ_2\in\mbR^{n\times(p-q)}$.
Then $\mathcal{C}(\bZ)=\mathcal{C}(\bX)$, and the linear model \eqref{eqn: by = bX bbeta + bepsilon} can be written as
\beq 
\by
=
\bX \bbeta + \bepsilon
=
\bZ  \begin{bmatrix}
    \bC \bbeta\\
    \bD \bbeta
\end{bmatrix}
+
\bepsilon
=
\bZ_1 \bC \bbeta + \bZ_2 \bD \bbeta + \bepsilon.
\eeq

Then the OLS estimator in the regression of $\by$ on
$\bZ$ satisfies
\beq
(\bZ\mt\bZ)^{-1}\bZ\mt\by
\overset{\eqref{eqn: transformed Z in regression}}{=}
\begin{bmatrix}
\bC\\
\bD
\end{bmatrix}
(\bX\mt\bX)^{-1}\bX\mt\by
=
\begin{pmatrix}
\bC\widehat{\bbeta}\\
\bD\widehat{\bbeta}
\end{pmatrix}.
\eeq
Let $\widetilde{\bZ}_1 = (\bI_n-\bP_{\bZ_2})\bZ_1$. Then by Proposition \ref{prop: FWL Thm} (i), we have
\beq
\bC\widehat{\bbeta}
=
(\widetilde{\bZ}_1\mt\widetilde{\bZ}_1)^{-1}
\widetilde{\bZ}_1\mt\by.
\eeq
Since
$\by=\bZ_1\bC\bbeta+\bZ_2\bD\bbeta+\bepsilon$,
$\widetilde{\bZ}_1\mt\bZ_2=\bzero$, and
$\widetilde{\bZ}_1\mt\bZ_1
=\widetilde{\bZ}_1\mt\widetilde{\bZ}_1$, we have
\beq
\bC(\widehat{\bbeta}-\bbeta)
=
(\widetilde{\bZ}_1\mt\widetilde{\bZ}_1)^{-1}
\widetilde{\bZ}_1\mt\bepsilon.
\eeq
Moreover, by Proposition \ref{prop: FWL Thm} (ii), the upper-left $q\times q$ block of $(\bZ\mt\bZ)^{-1}$ is $(\widetilde{\bZ}_1\mt\widetilde{\bZ}_1)^{-1}$.
On the other hand, by \eqref{eqn: transformed Z in regression}, we have
\beq
(\bZ\mt\bZ)^{-1}
=
\begin{bmatrix}
\bC\\
\bD
\end{bmatrix}
(\bX\mt\bX)^{-1}
\begin{bmatrix}
\bC\mt & \bD\mt
\end{bmatrix},
\eeq
whose upper-left block is
$\bC(\bX\mt\bX)^{-1}\bC\mt$.
Hence we have
\beq
\bC(\bX\mt\bX)^{-1}\bC\mt
=
(\widetilde{\bZ}_1\mt\widetilde{\bZ}_1)^{-1}.
\eeq
Thus,
\begin{eqnarray}
(\widehat{\bbeta}-\bbeta)\mt\bC\mt
\left\{
\bC(\bX\mt\bX)^{-1}\bC\mt
\right\}^{-1}
\bC(\widehat{\bbeta}-\bbeta)
=
\bepsilon\mt\bP_{\widetilde{\bZ}_1}\bepsilon
=
\bepsilon\mt
(\bP_{\bZ}-\bP_{\bZ_2})
\bepsilon,
\label{eqn: regression F numerator projection}
\end{eqnarray}
where the last equality follows from Proposition \ref{prop: FWL Thm} (iii).

Moreover, since $\mathcal{C}(\bX)=\mathcal{C}(\bZ)$, we have
\beq
\by\mt(\bI_n-\bP_{\bX})\by
=
\bepsilon\mt(\bI_n-\bP_{\bZ})\bepsilon.
\eeq
Thus the $F$ statistic in \eqref{eqn: F statistic in regression} can be written as
\begin{equation}
F
=
\frac{
\bepsilon\mt(\bP_{\bZ}-\bP_{\bZ_2})\bepsilon/q
}{
\bepsilon\mt(\bI_n-\bP_{\bZ})\bepsilon/(n-p)
}.
\label{eqn: F stat distn Z}
\end{equation}

Under the conditions of Theorem \ref{thm::classic-F-ols-2}, $\bX$, and therefore $\bZ$, is fixed, while $\bepsilon$ is spherically symmetric with $\p\{\bepsilon=\bzero_n\}=0$.
The result then follows immediately from
Theorem \ref{thm::classic-F-2}.

Under the conditions of Theorem \ref{thm::F-ols-new}, for any orthogonal matrix $\bPhi$, by \eqref{eqn: transformed Z in regression} and Condition \ref{(X2)}, we have $\bPhi \bZ \overset{d}{=} \bZ$, i.e., $\bZ$ is left-spherically symmetric. 
Also, $\rank(\bZ)=\rank(\bX)=p$ almost surely.
Since $\bepsilon \neq \bzero_n$ is fixed, the result follows immediately from Theorem \ref{thm::F-new}.

\section{Heuristics on Permutation and Spherical Symmetry}
\label{sec: permutation and spherical symmetry}

In this section, we use a geometric interpretation of permutation to extend our duality perspective of the $F$ statistic to settings where spherical symmetry is approximated by permutation symmetry.
This echoes Fisher's randomization argument for analysis of variance (ANOVA).
As discussed in \citet[p.~2]{box1980fisher}, Fisher viewed random allocation as creating a symmetry resembling the spherical symmetry supplied by normal theory, thereby motivating the approximate validity of the usual variance-ratio test.

To illustrate Fisher's intuition, we first provide a geometric interpretation of the permutation distribution.
For a fixed nonconstant vector, such as $\by=(1,2,\dots,n)\mt\in\mathbb{R}^n$, all permutations preserve, among other quantities, the coordinate sum $\sum_{i=1}^n y_i$ and the Euclidean norm $\|\by\|_2$.
The fixed sum $\sum_{i=1}^n y_i$ confines all permuted vectors to an $(n-1)$-dimensional affine hyperplane.
After centering, this hyperplane becomes the mean-zero subspace $\{\bu:\bone_n\mt\bu=0\}$.
Moreover, all centered permutation vectors have the same norm $\|\by-\bar{y}\bone_n\|_2$.
Thus they lie on a sphere within this $(n-1)$-dimensional subspace, and hence the sphere has dimension $n-2$.

Permutations preserve additional quantities, such as $\sum_{i=1}^n y_i^k$, $k \geq 3$, so the permutation vectors form only a finite subset of this sphere.
The spherical heuristic amounts to replacing this discrete set by a continuous distribution on the surrounding sphere.
This geometric representation of permutations is discussed by \citet[pp.~278--280]{mccullagh1993spheres}.

\begin{figure}[h!]
\centering
\includegraphics[width=0.95\textwidth]{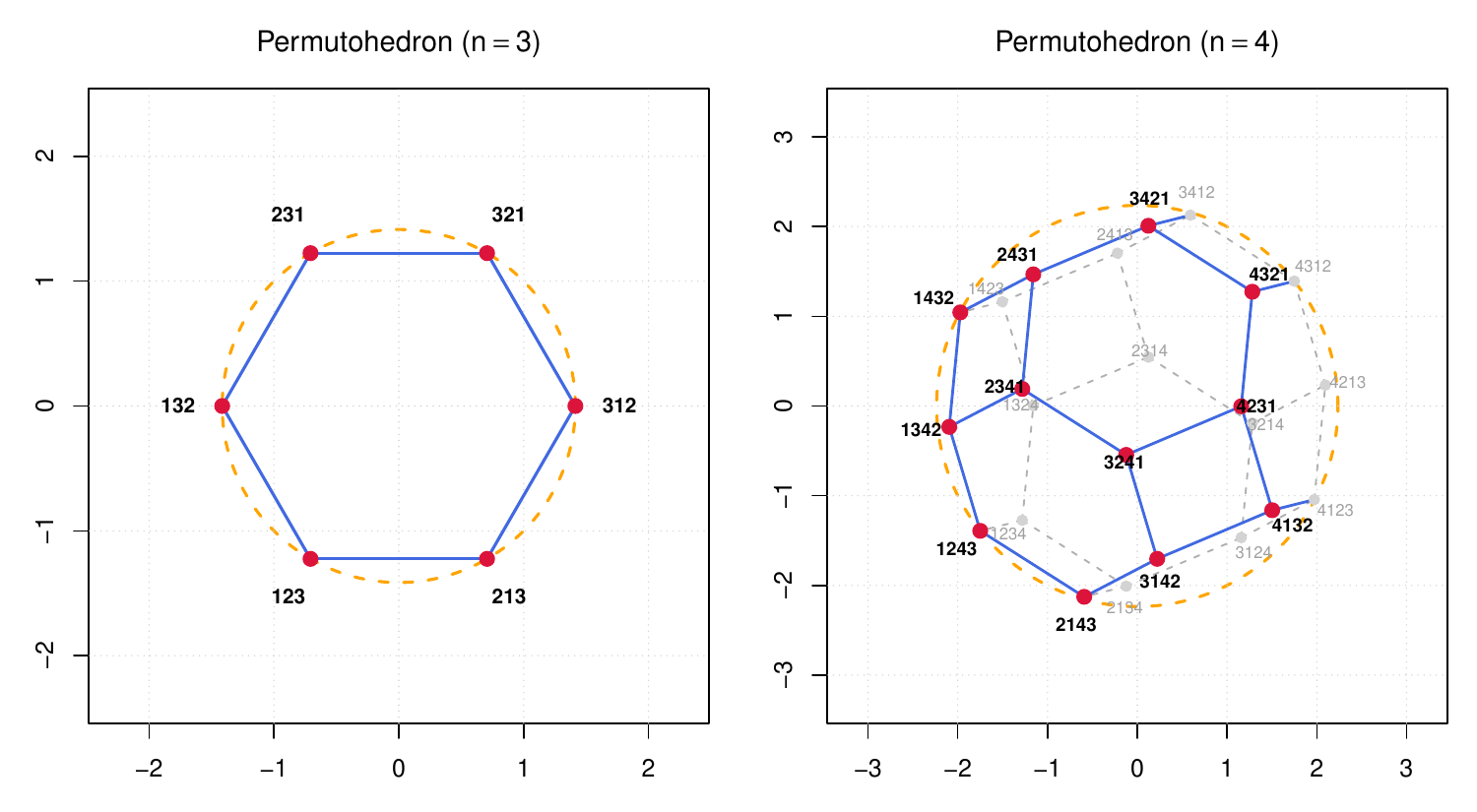}
\caption{Geometric representations of all permutations of
$(1,2,3)\mt$ (left) and $(1,2,3,4)\mt$ (right), after centering.
After rescaling, the six points for $n=3$ lie on $S^1$, whereas the 24 points for $n=4$ lie on $S^2$.
The figure follows the geometric representation in \citet[p.~279]{mccullagh1993spheres} and is reproduced with the help of AI.}
\label{fig: permutation sphere}
\end{figure}

The convex hull of the permutation vectors is known as a \emph{permutohedron}.
Figure~\ref{fig: permutation sphere} illustrates two simple examples.
When $\by=(1,2,3)\mt$, the six permutation vectors are the vertices of a regular hexagon, whose vertices lie on a circle.
When $\by=(1,2,3,4)\mt$, the 24 permutation vectors are the vertices of a three-dimensional permutohedron, whose vertices lie on a sphere.
Under a uniform random permutation, equal probability is assigned to these finitely many permutation vectors.

To connect this geometric view with the duality in our results, consider a randomized experiment with $K$ treatment groups of sizes $n_1,\ldots,n_K$, where $\sum_{k=1}^K n_k=n$.
Let $(G_1,\dots,G_n)\in \{1,\dots,K\}^n$ be the treatment group assignment vector,
and let 
\beq 
\bZ_0  = \begin{bmatrix}
    \mathcal{I}\{ G_1 = 1 \} & \dots & \mathcal{I}\{ G_1 = K \}\\
    \vdots & & \vdots\\
    \mathcal{I}\{ G_n = 1 \} & \dots & \mathcal{I}\{ G_n = K \}
\end{bmatrix} \in \mathbb{R}^{n \times K}
\eeq
be the fixed group-indicator matrix.
Let $\bQ \in \mathbb{R}^{n \times n}$ be a uniformly distributed random permutation matrix, so that a randomized assignment can be represented as $\bZ=\bQ\bZ_0$.
The usual one-way ANOVA statistic can be written as
\begin{equation}
F(\by,\bZ)
=
\frac{
\by\mt(\bP_{\bZ}-\bP_{\bone_n})\by/(K-1)
}{
\by\mt(\bI_n-\bP_{\bZ})\by/(n-K)
},
\qquad
\bP_{\bone_n}
=
\frac{\bone_n\bone_n\mt}{n}.
\label{eqn: permutation ANOVA F}
\end{equation}
Because $\bP_{\bQ\bZ_0}=\bQ\bP_{\bZ_0}\bQ\mt$ and $\bQ\bone_n=\bone_n$, we have the exact identity
\begin{equation}
F(\by,\bQ\bZ_0)
=
F(\bQ\mt\by,\bZ_0).
\label{eqn: permutation ANOVA duality}
\end{equation}

Equation~\eqref{eqn: permutation ANOVA duality} shows two equivalent views of the same randomization.
On its right-hand side, the treatment groups are fixed and the outcomes are randomly permuted. Heuristically, replacing the finite set of permuted outcome vectors by a uniform distribution on the surrounding sphere suggests the classical spherical-symmetry view underlying Theorem~\ref{thm::classic-F-2}.
On its left-hand side, the outcomes are fixed and the treatment assignment is randomly permuted. Heuristically, replacing the discrete row permutations of $\bZ_0$ by continuous symmetric rotations suggests the dual random-design view underlying Theorem~\ref{thm::F-new}.
Thus the two equivalent views of randomization echo our duality of the $F$ distribution: the randomness can be viewed as acting either on the vector being projected or on the space onto which it is projected.

As pointed out by \citet{stigler2004discussion, stigler2016seven},
Fisher’s initial justification of the $F$ distribution for ANOVA was that, with sufficiently rich randomization, discrete symmetry can be viewed as an approximation to continuous spherical symmetry.
The rigorous randomization-based theory relies on limiting theorems of the permutation statistics \citep{ding2018randomization}.

\section{Discussion}

We have reviewed the classic result on the $F$ distribution based on spherically symmetric random errors and presented a dual result based on a left-spherically symmetric random matrix. We also extend the result to linear regression. Our results complement recent discussion of statistical inference issues in fixed design and random design linear regression \citep{abadie2020sampling, zhang2025random}. In Appendix~\ref{sec: F test with mixed regressors}, we extend the result to mixed designs with both fixed and random regressors \citep{chetverikov2023standard, gao2026misspecified}.



\clearpage
\section*{Supplementary Materials}
\setcounter{section}{0}
\renewcommand\thesection{\Alph{section}}
\renewcommand\thesubsection{\thesection.\arabic{subsection}}


Section \ref{sec: F test with mixed regressors} presents the duality of the $F$ statistic with mixed regressors in linear regression.
Section \ref{sec: proof of sec F test with mixed regressors} presents the proofs for the results in Section \ref{sec: F test with mixed regressors}.



\setcounter{equation}{0}
 \makeatletter
 \renewcommand{\theequation}{A.\@arabic\c@equation}
\makeatother

\setcounter{lemma}{0}
 \makeatletter
\renewcommand{\thelemma}{A.\@arabic\c@lemma}
\makeatother

\section{Extension: Mixed Regressors}
\label{sec: F test with mixed regressors}

This section investigates the validity of the $F$ statistic in settings where $\bX$ consists of both fixed and random regressors.


Specifically, in the linear model \eqref{eqn: by = bX bbeta + bepsilon},
let $\bX = [\bX_0\ \bX_1] \in \mathbb{R}^{n \times p}$ with $\mathrm{rank}(\bX) = p$, where $\bX_0 \in \mathbb{R}^{n \times s}$ is fixed and $\bX_1 \in \mathbb{R}^{n \times d}$ is random with $p = s + d$.
Let $\bbeta = [\bbeta_0\mt\ \bbeta_1\mt]\mt \in \mathbb{R}^{p}$ where $\bbeta_0 \in \mathbb{R}^s$ and $\bbeta_1 \in \mathbb{R}^d$ correspond to $\bX_0$ and $\bX_1$, respectively.
So the model in \eqref{eqn: by = bX bbeta + bepsilon} becomes
\begin{equation}
\by = \bX_0 \bbeta_0 + \bX_1 \bbeta_1 + \bepsilon.
\label{eqn: by = bX0beta0 + bX1beta1 + beps}    
\end{equation}
The OLS estimator in \eqref{eqn: OLS estimator} can be written as $\widehat{\bbeta} = [\widehat{\bbeta}_0\mt\ \widehat{\bbeta}_1\mt]\mt$.

Recall that the general form of $F$ statistic in \eqref{eqn: F statistic in regression} is designed for $\bC \bbeta$.
In this section, our focus is on inference related to $\bbeta_1$, the coefficients of random regressors.
Accordingly, the general matrix $\bC$ takes the form
\begin{equation}
\bC = [ \mathbf{O}_{q \times s}, \bC_1 ]\in \mathbb{R}^{q \times (s+d) },
\label{eqn: bC = (O bC1)}
\end{equation} 
where  $\bC_1 \in \mathbb{R}^{q \times d}$ is any fixed matrix with $\mathrm{rank}(\bC_1) = q$.
In this case, we have $\bC \widehat{\bbeta} = \bC_1 \widehat{\bbeta}_1$, $\bC \bbeta = \bC_1 \bbeta_1$, and the $F$ statistic in \eqref{eqn: F statistic in regression} becomes
\begin{equation}
F =
\frac{ \{ (\widehat{\bbeta}_1-\bbeta_1)\mt \bC_1\mt \lmk \bC (\bX\mt\bX)\inverse \bC\mt \rmk\inverse \bC_1 ( \widehat{\bbeta}_1 - \bbeta_1 ) \} / q }{  [ \by\mt \lsk \bI_n-\bP_{\bX} \rsk \by  ] / (n-p)}.
\label{eqn: F statistic with C1}
\end{equation}
Suppose $1 \leq q \leq d \leq d+s = p < n$.
If $s = 0$, $\bX$ is reduced to $\bX_1$ and $\bC = \bC_1$, so that $F$ in \eqref{eqn: F statistic with C1} reduces to \eqref{eqn: F statistic in regression}.

Let $\boldsymbol{\Gamma} \in \mathbb{R}^{n \times (n - s)}$ be a matrix whose columns form an orthonormal basis of $\mathcal{C}(\bX_0)^{\perp}$.
We define the projected quantities:
\begin{equation}
\by^\star = \bGamma\mt \by \in \mathbb{R}^{n-s}, 
\ 
\bX_1^\star = \bGamma\mt \bX_1  \in \mathbb{R}^{(n-s) \times d},
\  
\bepsilon^\star = \bGamma\mt \bepsilon \in \mathbb{R}^{n-s}.
\label{eqn: bystar bXstar bepsstar}
\end{equation}
Then $F$ statistic in \eqref{eqn: F statistic with C1} takes a simplified form, as shown in the following algebraic result:

\begin{proposition}
Under \eqref{eqn: by = bX0beta0 + bX1beta1 + beps} with $\mathrm{rank}(\bX) = s + d$, let $\bGamma \in \mathbb{R}^{n \times (n-s)}$ be formed from any orthonormal basis of $\mathcal{C}(\bX_0)^{\perp}$. Then
\begin{enumerate}
    \item[(i)] $\by^\star = \bX_1^\star \bbeta_1 + \bepsilon^\star$;
    \item[(ii)] Define $\widehat{\bbeta}_1^\star = ( \bX_1^{\star \T} \bX_1^\star )\inverse \bX_1^{\star \T} \by^\star$.
    Then the OLS estimator under \eqref{eqn: by = bX0beta0 + bX1beta1 + beps} satisfies $\widehat{\bbeta}_1 = \widehat{\bbeta}_1^\star$;
    \item[(iii)] The $F$ statistic in \eqref{eqn: F statistic with C1} satisfies
\end{enumerate}
\begin{equation}
F 
=
\frac{  \{ (\widehat{\bbeta}_1^\star -\bbeta_1)\mt \bC_1\mt [ \bC_1 (\bX_1^{\star\T}\bX_1^\star)\inverse \bC_1\mt ]\inverse \bC_1 ( \widehat{\bbeta}_1^\star-\bbeta_1) \} / q  }{   \{ \by^{\star\T} \lsk \bI_{n-s}-\bP_{\bX_1^\star} \} \by^\star  \rbk / (n-p)  }.
\label{eqn: F statistic X1star ystar}    
\end{equation}
\label{prop: Gammat model}
\end{proposition}

The space $\mathcal{C}(\bX_0)^{\perp}$ and its orthonormal basis $\bGamma$ represent the information not captured by the fixed regressors in model~\eqref{eqn: by = bX0beta0 + bX1beta1 + beps}.
Thus,  $\by^\star$, $\bX_1^\star$, and $\bepsilon^\star$ contain the new information not explained by the fixed regressors.
Proposition \ref{prop: Gammat model} shows that the inference on $\bbeta_1$ under the full model \eqref{eqn: by = bX0beta0 + bX1beta1 + beps} can be equivalently conducted under the reduced model $\by^\star = \bX_1^\star \bbeta_1 + \bepsilon^\star$.
By Theorems \ref{thm::classic-F-ols-2} and \ref{thm::F-ols-new}, this result highlights a duality between $\bX_1^\star$ and $\bepsilon^\star$ in the inference of $\bbeta_1$.

\begin{corollary}
[Duality of $F$ statistic with mixed regressors]
Under model \eqref{eqn: by = bX0beta0 + bX1beta1 + beps}, 
suppose $\bepsilon \Perp \bX_1$, $\p \{ \mathrm{rank}( [\bX_0, \bX_1] ) < s+d \} = 0$, and $\p\{ \bepsilon \in \mathcal{C}(\bX_0) \} = 0$.
Assume that $\bC$ in \eqref{eqn: bC = (O bC1)} satisfies $\mathrm{rank}(\bC_1) = q$ and $1 \leq q \leq d \leq d+s = p < n$.
\begin{itemize}
    \item[(a)] If there exists $\bGamma \in \mathbb{R}^{n \times (n - s)}$ formed by an orthonormal basis of $\mathcal{C} (\bX_0)^{\perp}$ such that $\bepsilon^\star$ in \eqref{eqn: bystar bXstar bepsstar} is a spherically symmetric random vector as defined in \ref{(W1-1)}, then the $F$ statistic in \eqref{eqn: F statistic with C1} follows $F(q, n-p)$.
    \item[(b)] If there exists $\bGamma \in \mathbb{R}^{n \times (n - s)}$ formed by an orthonormal basis of $\mathcal{C} (\bX_0)^{\perp}$ such that $\bX_1^\star$ in \eqref{eqn: bystar bXstar bepsstar} is a left-spherically symmetric random matrix as defined in \ref{(Z2)}, then the $F$ statistic in \eqref{eqn: F statistic with C1} follows $F(q, n-p)$.
\end{itemize}
\label{coro: F statistic mixed X}
\end{corollary}

\setcounter{equation}{0}
 \makeatletter
 \renewcommand{\theequation}{B.\@arabic\c@equation}
\makeatother

\setcounter{lemma}{0}
 \makeatletter
\renewcommand{\thelemma}{B.\@arabic\c@lemma}
\makeatother

\section{Additional Proofs and Propositions}
\label{sec: proof of sec F test with mixed regressors}

\subsection{Proof of Proposition \ref{prop: Gammat model}}

By the definition of $\bGamma \in \mathbb{R}^{n \times (n-s)}$, we have 
$
\bI_n - \bP_{\bX_0} = \bGamma \bGamma\mt  
$.
Multiplying $\bGamma\mt$ on both sides of \eqref{eqn: by = bX0beta0 + bX1beta1 + beps}, we have
$
\by^\star = \bX_1^\star \bbeta_1 + \bepsilon^\star
$,
which completes the proof of Part (i) of Proposition \ref{prop: Gammat model}.

Based on 
$
\by^\star = \bX_1^\star \bbeta_1 + \bepsilon^\star
$, the OLS estimator is
$
\widehat{\bbeta}_1^\star 
=
(\bX_1^{\star \T} \bX_1^\star)\inverse \bX_1^{\star \T} \by^\star
$.
By Proposition \ref{prop: FWL Thm} (i) and $\bGamma\bGamma\mt = \bI_n - \bP_{\bX_0}$, we have
$
\widehat{\bbeta}_1
=
[\bX_1\mt (\bI_n - \bP_{\bX_0}) \bX_1]\inverse \bX_1\mt (\bI_n - \bP_{\bX_0}) \by = \widehat{\bbeta}_1^\star
$,
which completes the proof of Part (ii) of Proposition \ref{prop: Gammat model}.

The fact that $\widehat{\bbeta}_1 = \widehat{\bbeta}_1^\star$ also implies that
\begin{equation}
\bC (\widehat{\bbeta} - \bbeta)
=
\bC_1 (\widehat{\bbeta}_1 - \bbeta_1)
=
\bC_1 (\widehat{\bbeta}_1^\star - \bbeta_1).
\label{eqn: C (betahat-beta) = C1 (betahat1star - beta1)}    
\end{equation}

From Proposition \ref{prop: FWL Thm} (ii) and the definition of $\bC$ in \eqref{eqn: bC = (O bC1)}, we have
\begin{eqnarray}
\bC (\bX\mt\bX)\inverse \bC\mt 
=
\bC_1 [ \bX_1\mt (\bI_n - \bP_{\bX_0}) \bX_1 ]\inverse \bC_1\mt 
=
\bC_1 (\bX_1^{\star \T} \bX_1^\star)\inverse \bC_1\mt,
\label{eqn: C solve(XtX) Ct = C1 solve(X1start X1star) C1t}    
\end{eqnarray}
where the second equality follows from $\bGamma\bGamma\mt = \bI_n - \bP_{\bX_0}$.

Since $\bP_{(\bI_n - \bP_{\bX_0}) \bX_1}$ and $\bP_{\bX_1^\star}$ are the projection matrices onto $\mathcal{C} ( (\bI_n - \bP_{\bX_0}) \bX_1 )$ and $\mathcal{C} (\bX_1^\star)$, respectively, we have
\begin{eqnarray}
\bP_{(\bI_n - \bP_{\bX_0}) \bX_1}
=
(\bI_n - \bP_{\bX_0}) \bX_1 [ \bX_1\mt (\bI_n - \bP_{\bX_0}) \bX_1 ]^{-1} \bX_1\mt (\bI_n - \bP_{\bX_0})
=
\bGamma \bP_{\bX_1^\star} \bGamma\mt.
\label{eqn: P(In - P(X0)X1) = Gamma P(X1star) Gammat}  
\end{eqnarray}
Then by Proposition \ref{prop: FWL Thm} (iv), we have
\begin{equation}
(\bI_n - \bP_{\bX}) \by
=
[\bI_n - \bP_{(\bI_n - \bP_{\bX_0}) \bX_1}] (\bI_n - \bP_{\bX_0}) \by
=
(\bI_n - \bGamma \bP_{\bX_1^\star} \bGamma\mt) \bGamma \by^\star 
= 
\bGamma (\bI_{n-s} - \bP_{\bX_1^\star}) \by^\star,
\label{eqn: (I - P(X))y = Gamma (I-P(X1star)) ystar}
\end{equation}
where the second equality follows from \eqref{eqn: P(In - P(X0)X1) = Gamma P(X1star) Gammat} and $\by^\star = \bGamma\mt \by$.

Combining \eqref{eqn: C (betahat-beta) = C1 (betahat1star - beta1)}, \eqref{eqn: C solve(XtX) Ct = C1 solve(X1start X1star) C1t}, and \eqref{eqn: (I - P(X))y = Gamma (I-P(X1star)) ystar} completes the proof of Part (iii) of Proposition \ref{prop: Gammat model}.
\qed

\subsection{Proof of Corollary \ref{coro: F statistic mixed X}}

We apply Theorems \ref{thm::classic-F-ols-2} and \ref{thm::F-ols-new} based on the $F$ statistic in \eqref{eqn: F statistic X1star ystar} and $
\by^\star = \bX_1^\star \bbeta_1 + \bepsilon^\star
$.
By Proposition \ref{prop: FWL Thm} (i) and $\mathrm{rank}([\bX_0,\bX_1]) = p$ almost surely, $(\bI_n - \bP_{\bX_0}) \bX_1 = \bGamma \bGamma\mt \bX_1 = \bGamma \bX_1^\star$ is of full rank almost surely, which implies that $\p \{ \mathrm{rank}(\bX_1^\star) < d \} = 0$.
By the Orthogonal Projection Theorem, $\bepsilon \notin \mathcal{C}(\bX_0)$ is equivalent to $\bGamma\mt \bepsilon \neq \bzero_{n-s}$, $\ie$, $\bepsilon^\star \neq \bzero_{n-s}$.
Thus $\p \{ \bepsilon^\star = \bzero_{n-s} \} = \p \{ \bepsilon \in \mathcal{C}(\bX_0) \} = 0$.

For Part (a) of Corollary \ref{coro: F statistic mixed X},  by $\bepsilon^\star \Perp \bX_1^\star$,  Theorem \ref{thm::classic-F-ols-2} implies that $(F \mid \bX_1^\star) \sim F(q, n-p)$ almost surely.
Then the marginal distribution of the $F$ statistic in \eqref{eqn: F statistic X1star ystar} follows $F(q, n-p)$.
For Part (b) of Corollary \ref{coro: F statistic mixed X}, Theorem \ref{thm::F-ols-new} implies that $(F \mid \bepsilon^\star) \sim F(q, n-p)$ almost surely.
Then the marginal distribution of the $F$ statistic in \eqref{eqn: F statistic X1star ystar} also follows $F(q, n-p)$.
\qed


\section*{Acknowledgments}

Professors Don Rubin, Steve Stigler, and Peter McCullagh gave us helpful historical background and references.  
The main results of this paper appeared in \citet[][Chapter 3]{Zhang2025thesis}. When we were polishing this paper, we discovered overlapping results in \citet[][Lemma 1]{wen2025residual} through \citet[][Proposition 1]{xia2026regression}. Our presentation is more pedagogical, and we have additional interests in the heuristics on permutation and spherical symmetry and the theory with mixed regressors.

\bibliography{bibliography.bib}

\end{document}